\documentclass[aos,MSNbibl,nameyear,dvips]{arximspdf}
\usepackage{dcolumn}
\usepackage{graphicx}

\doi{10.1214/14-AOS1207} 
\volume{42}
\issue{2}
\pubyear{2014}
\firstpage{789}
\lastpage{817}

\makeatletter
\def\mid{|}
\newcolumntype{d}[1]{D{.}{.}{#1}}
\newcommand{\rrVert}{\Vert}
\newcommand{\llVert}{\Vert}
\renewcommand{\citep}[1]{(\citeauthor{#1} \citeyear{#1})}
\newcommand{\eqref}[1]{(\ref{#1})}
\newcommand{\cprob}{\stackrel{\mathrm{P}}{\longrightarrow}}
\newcommand{\eqdist}{\stackrel{ d}{=}}

\newtheorem{thmm}{Theorem}[section]

\newtheorem{lem}{Lemma}[section]
\newtheorem{cor}{Corollary}[section]
\newproclaim{definition}[thmm]{Definition}
\newproclaim{notation}[thmm]{Notation}
\newproclaim{example}[thmm]{Example}
\newproclaim{condition}{Condition}[section]
\newproclaim{remark}{Remark}
\makeatother

\begin{document}
\begin{frontmatter}

\title{Bayesian variable selection with shrinking and diffusing
priors\thanksref{T1}}
\runtitle{Variable selection with shrinking and diffusing priors}
\thankstext{T1}{Supported in part by the NSF Grants DMS-12-37234 and
DMS-13-07566, and by the National Natural Science Foundation of China Grant
11129101.}

\begin{aug}
\author{\fnms{Naveen Naidu}~\snm{Narisetty}\ead[label=e1]{naveennn@umich.edu}}
\and
\author{\fnms{Xuming}~\snm{He}\corref{}\ead[label=e2]{xmhe@umich.edu}}
\runauthor{N. N. Narisetty and X. He}
\affiliation{University of Michigan}

\address{Department of Statistics\\
University of Michigan\\
Ann Arbor, Michigan 48109\\
USA\\
\printead{e1}\\
\phantom{E-mail:\ }\printead*{e2}}
\end{aug}

\received{\smonth{5} \syear{2013}}
\revised{\smonth{1} \syear{2014}}

%
\begin{abstract}
We consider a Bayesian approach to variable selection in the presence
of high dimensional covariates based on a hierarchical model that
places prior distributions on the regression coefficients as well as on
the model space. We adopt the well-known spike and slab Gaussian priors
with a distinct feature, that is, the prior variances depend on the
sample size through which appropriate shrinkage can be achieved. We
show the strong selection consistency of the proposed method in the
sense that the posterior probability of the true model converges to one
even when the number of covariates grows nearly exponentially with the
sample size. This is arguably the strongest selection consistency
result that has been available in the Bayesian variable selection
literature; yet the proposed method can be carried out through
posterior sampling with a simple Gibbs sampler. Furthermore, we argue
that the proposed method is asymptotically similar to model selection
with the $L_0$ penalty. We also demonstrate through empirical work the
fine performance of the proposed approach relative to some state of the
art alternatives.
\end{abstract}

%


\end{frontmatter}

\section{Introduction}\label{intro-sec}

We consider the linear regression setup with high dimensional
covariates where the number of covariates $p$ can be large relative to
the sample size $n$. When $p > n$, the estimation problem is ill-posed
without performing variable selection. A natural assumption to limit
the number of parameters in high dimensional settings is that the
regression function (i.e., the conditional mean) is sparse in the sense
that only a small number of covariates (called active covariates) have
nonzero coefficients. We aim to develop a new Bayesian methodology for
selecting the active covariates that is asymptotically consistent and
computationally convenient. A~large number of methods have been
proposed for variable selection in the literature from both frequentist
and Bayesian viewpoints. Many frequentist methods based on penalization
have been proposed following the well-known least absolute shrinkage
and selection operator [LASSO, \citet{Tibs96}]. We mention the smoothly
clipped absolute deviation [SCAD, \citet{Fan01}], adaptive LASSO [\citet
{Zou06}], octagonal shrinkage and clustering algorithm for regression
[OSCAR, \citet{Bondell08}] and the Dantzig selector [\citet{Candes07};
\citet{James09}] just to name a few. \citet{Fan10} provided a selective
overview of high dimensional variable selection methods. Various
authors reported inconsistency of LASSO and its poor performance for
variable selection under high dimensional settings; see \citet{Zou06}
and \citet{Johnson12}. On the other hand, several penalization based
methods were shown to have the oracle property [\citet{Fan01}] under
some restrictions on~$p$. For example, \citet{Fan04} and \citet{Huang07}
showed the oracle property for some nonconcave penalized likelihood
methods when $p = O(n^{1/3})$ and $p = o(n)$, respectively. \citet
{Shen12} showed that $L_0$ penalized likelihood method has the oracle
property under exponentially large $p = e^{o(n)}$.

Many Bayesian methods have also been proposed for variable selection
including the stochastic search variable selection [\citet{George93}],
empirical Bayes variable selection [\citet{George00}], spike and slab
selection method [\citet{Ishwaran05}], penalized credible regions
[\citet
{Bondell12}], nonlocal prior method [\citet{Johnson12}], among others.
We shall describe the typical framework used for Bayesian variable
selection methods before discussing their theoretical properties.

We use the standard notation $Y_{n \times1} = X_{n \times p} \beta_{p
\times1} + \varepsilon_{n \times1}$ to represent the linear regression
model. Bayesian variable selection methods usually introduce latent
binary variables for each of the covariates to be denoted by $Z = (Z_1,
\ldots, Z_p)$. The idea is that each $Z_i$ would indicate whether the
$i$th covariate is active in the model or not. For this reason, the
prior distribution on the regression coefficient $\beta_i$ under $Z_i =
0$ is usually a point mass at zero, but a diffused (noninformative)
prior under $Z_i =1$. The concentrated prior of $\beta_i$ under $Z_i =
0$ is referred to as the spike prior, and the diffused prior under $Z_i
= 1$ is called the slab prior. Further, a prior distribution on the
binary random vector $Z$ is assumed, which can be interpreted as a
prior distribution on the space of models. A Bayesian variable
selection method then selects the model with the highest posterior
probability. Various selection procedures with this structure have been
proposed; they essentially differ in the form of the spike and slab
priors, or in the form of the prior on the model space.

\citet{Mitchell88} considered a uniform distribution for the slab prior.
\citet{George93} used the Gaussian distribution with a zero mean and a
small but fixed variance as the spike prior, and another Gaussian
distribution with a large variance as the slab prior. This allowed the
use of a Gibbs sampler to explore the posterior distribution of $Z$.
However, as we argue in Section~\ref{ortho-sec}, this prior
specification does not guarantee model selection consistency at any
fixed prior.
\citet{Ishwaran05} also used Gaussian spike and slab priors, but with
continuous bimodal priors for the variance of $\beta$ to alleviate the
difficulty of choosing specific prior parameters. More recently, \citet
{Ishwaran11} established the oracle property for the posterior mean as
$n$ converges to infinity (but $p$ is fixed) under certain conditions
on the prior variances. They noted that in the orthogonal design case,
a uniform complexity prior leads to correct complexity recovery (i.e.,
the expected size of the posterior model size converges to the true
model size) under weaker conditions on the prior variances.
In another development, \citet{Yang12} used shrinking priors to explore
commonality across quantiles in the context of Bayesian quantile
regression, but the use of such priors for achieving model selection
consistency has not been explored. In this paper, we continue to work
with the framework where both the spike and slab priors are Gaussian,
but our prior parameters depend explicitly on the sample size through
which appropriate shrinkage is achieved. We shall establish model
selection consistency properties for general design matrices while
allowing $p$ to grow with $n$ at a nearly exponential rate. In
particular, the strong selection consistency property we establish is a
stronger result for model selection than complexity recovery.

One of the most commonly used priors on the model space is the
independent prior given by $P[Z = z] = \prod_{i =1}^p w_i^{z_i} (1 -
w_i)^{z_i}$, where the marginal probabilities $w_i$ are usually taken
to be the same constant. However, when $p$ is diverging, this implies
that the prior probability on models with sizes of order less than $p$
goes to zero, which is against model sparsity. We consider marginal
probabilities $w_i$ in the order of $p^{-1}$, which will impose
vanishing prior probability on models of diverging size. \citet{Yuan05}
used a prior that depends on the Gram matrix to penalize models with
unnecessary covariates at the prior level. The vanishing prior
probability in our case achieves similar prior penalization.

A common notion of consistency for Bayesian variable selection is
defined in terms of pairwise Bayes factors, that is, the Bayes factor
of any under- or over-fitted model with respect to the true model goes
to zero. \citet{Moreno10} proved that intrinsic priors give pairwise
consistency when $p = O(n)$, and similar consistency of the Bayesian
information criterion [BIC, \citet{Schwarz78}] when $p = O(n^{\alpha}),
\alpha<1$. Another notion of consistency for both frequentist and
Bayesian methods is that the selected model equals the true model with
probability converging to one. We refer to this as selection
consistency. \citet{Bondell12} proposed a method based on penalized
credible regions that is shown to be selection consistent when $\log p
= O(n^c), c<1$. \citet{Johnson12} proposed a stronger consistency for
Bayesian methods under which the posterior probability of the true
model converges to one, which we shall refer to as strong selection
consistency. The authors used nonlocal distributions (distributions
with small probability mass close to zero) as slab priors, and proved
strong selection consistency when $p < n$. However, apart from the
limitation $p < n$, their method involves approximations of the
posterior distributions and an application of MCMC methods, which are
computationally intensive if at all feasible for modest size problems.

We make the following contributions to variable selection in this
article. We introduce shrinking and diffusing priors as spike and slab
priors, and establish strong selection consistency of the approach for
$p = e^{o(n)}$. This approach is computationally advantageous because a
standard Gibbs sampler can be used to sample from the posterior. In
addition, we find that the resultant selection on the model space is
closely related to the $L_0$ penalized likelihood function. The merits
of the $L_0$ penalty for variable selection have been discussed by many
authors including \citet{Schwarz78}, \citet{Liu07}, \citet{Dicker13},
\citet
{Kim12} and \citet{Shen12}.

We now outline the remaining sections of the paper as follows. The
first part of Section~\ref{model-sec} describes the model, conditions
on the prior parameters and motivation for these conditions. The latter
part describes our proposed methodology for variable selection based on
the proposed model. Section~\ref{ortho-sec} motivates the use of sample
size dependent prior parameters by considering orthogonal design
matrices, and provides insight into the variable selection mechanism
using those priors. Section~\ref{results-sec} presents our main results
on the convergence of the posterior distribution of the latent vector
$Z$, and the strong selection consistency of our model selection
methodology. Section~\ref{penalty-sec} provides an asymptotic
connection between the proposed method and the $L_0$ penalization.
Section~\ref{cond-sec} provides a discussion on the conditions assumed
for proving the results of Section~\ref{results-sec}. Some
computational aspects of the proposed method are noted in Section~\ref{compu-sec}. We present simulation studies in Section~\ref{simul-sec}
to illustrate how the proposed method compares with some existing
methods. Application to a gene expression data set is given in
Section~\ref{data-sec}, followed by a conclusion in Section~\ref{conclusion-sec}. Section~\ref{proofs-sec} provides proofs of some
results not given in the earlier sections.

\section{The model} \label{model-sec}

From now on, we use $p_n$ to denote the number of covariates to
indicate that it grows with $n$. Consider the $n \times1$ response
vector $Y$, and the $n\times p_n$ design matrix $X$ corresponding to
the $p_n$ covariates of interest. Let $\beta$ be the regression vector,
that is, the conditional mean of $Y$ given $X$ is given by $X\beta$. We
assume that $\beta$ is sparse in the sense that only a few components
of $\beta$ are nonzero; this sparsity assumption can be relaxed as in
Condition~\ref{model-cond}. Our goal is to identify the nonzero
coefficients to learn about the active covariates. We describe our
working model as follows:
%
\begin{eqnarray}\label{modeleq}
\label{model} &&\hspace*{12pt}Y \mid\bigl( X, \beta, \sigma^2 \bigr) \sim N\bigl(X
\beta, \sigma^2 I\bigr),
\nonumber
\\
&&\beta_i \mid \bigl( \sigma^2, Z_i = 0
\bigr) \sim N\bigl(0, \sigma^2 {\tau ^2_{0,n}}
\bigr),
\nonumber
\\
&&\beta_i \mid \bigl(\sigma^2, Z_i = 1
\bigr) \sim N\bigl(0,\sigma^2 {\tau ^2_{1,n}}
\bigr),
\\
&&P(Z_i = 1 ) =1 - P(Z_i = 0 ) = q_n,
\nonumber
\\
&&\sigma^2 \sim \operatorname{IG}(\alpha_1,\alpha_2),\nonumber
\end{eqnarray}
where $i$ runs from 1 to $p_n$, $q_n, \tau_{0,n}, \tau_{1,n}$ are
constants that depend on $n$, and $\operatorname{IG}(\alpha_1,\alpha_2)$ is the
Inverse Gamma distribution with shape parameter $\alpha_1$ and scale
parameter $\alpha_2$.

The intuition behind this set-up is that the covariates with zero or
very small coefficients will be identified with zero $Z$ values, and the
active covariates will be classified as $Z = 1$. We use the posterior
probabilities of the latent variables $Z$ to identify the active
covariates.

\textit{Notation}:
We now introduce the following notation to be used throughout the
paper.

\textit{Rates}: For sequences $a_n$ and $b_n$, $a_n \sim b_n$ means
$\frac{a_n}{b_n} \rightarrow c$ for some constant $c >0$, $a_n \succeq
b_n$ (or $b_n \preceq a_n$) means $b_n = O(a_n)$, and $a_n \succ b_n$
(or $b_n \prec a_n$) means $b_n = o(a_n)$.

\textit{Convergence}: Convergence in probability is denoted by $\cprob
$, and equivalence in distribution is denoted by $\eqdist$.

\textit{Models}: We use $k$ to index an arbitrary model which is viewed
as a $p_n \times1$ binary vector. The $i$th entry $k_i$ of $k$
indicates whether the $i$th covariate is active (1) or not (0). We
use $X_k$ as the design matrix corresponding to the model $k$, and
$\beta_k$ to denote the corresponding regression coefficients. In
addition, $t$ is used to represent the true model.

\textit{Model operations}: We use $|k|$ to represent the size of the
model $k$. For two models $k$ and $j$, the operations $k \vee j$ and $k
\wedge j$ denote entry-wise maximum and minimum, respectively.
Similarly, $k^c = \mathbf{1} - k$ is entrywise operation, where
$\mathbf
{1}$ is the vector of 1's. We also use the notation $k \supset j$ (or
$k \geq j$) to denote that the model $k$ includes all the covariates in
model $j$, and $k \not\supset j$ otherwise.

\textit{Eigenvalues}: We use $\phi_{\mathrm{min}}(A)$ and $\phi_{\mathrm{max}}(A)$ to
denote the minimum and maximum eigenvalues, respectively, and $\phi
_{\mathrm{min}}^{\#}(A)$ to denote the minimum nonzero eigenvalue (MNEV) of the
matrix $A$. Moreover, we use $\lambda_M^n$ to be the maximum eigenvalue
of the Gram matrix $X'X/n$, and for $\nu> 0$, we define
\[
m_n(\nu) = p_n \wedge
\frac{n}{(2 + \nu)\log p_n }\quad \mbox{and}\quad \lambda_m^n (\nu): = \inf
_ {|k| \leq m_n(\nu)} \phi_{\mathrm{min}}^{\#} \biggl(
\frac{X_k'
X_k}{n} \biggr).
\]

\textit{Matrix inequalities}: For square matrices $A$ and $B$ of the
same order, $A \geq B$ or $ (A - B) \geq0$ means that $(A-B)$ is
positive semidefinite.

\textit{Residual sum of squares}: We define $\tilde{R}_k = Y' (I -
X(D_k + X'X)^{-1}X') Y$, where $D_k = \operatorname{Diag} (k \tau_{1n}^{-2} +
(\mathbf
{1} -k) \tau_{0n}^{-2})$. $\tilde{R}_k$ approximates the usual residual
sum of squares $R_k^* = Y' ( I - P_k  )Y$, where $P_k$ is the
projection matrix corresponding to the model~$k$.

\textit{Generic constants}: We use $c'$ and $w'$ to denote generic
positive constants that can take different values each time they appear.

\subsection{Prior parameters} \label{prior-sec} We consider ${\tau
^2_{0,n}} \rightarrow0$ and ${\tau^2_{1,n}} \rightarrow\infty$ as $n$
goes to $\infty$, where the rates of convergence depend on $n$ and
$p_n$. To be specific, we assume that for some $\nu>0$, and $\delta>0$,
\[
n \tau_{0n}^2
\lambda_M^n= o(1)\quad \mbox{and}\quad n \tau_{1n}^{2}
\lambda _m^n (\nu) \sim \bigl(n \vee
p_n^{2 + 2 \delta} \bigr).
\]
As will be seen later, these rates ensure desired model selection
consistency for any $\delta> 0$, where larger values of $\delta$ will
correspond to higher penalization and vice versa.

Note that the variance $\tau_{0n}^{2}$ depends on the sample size $n$
and the scale of the Gram matrix. Since the prior distribution of a
coefficient under $Z=0$ is mostly concentrated in
\[
\biggl(-\frac{3 \sigma}{\sqrt{n \lambda_M^n }}, \frac{3 \sigma
}{\sqrt{n
\lambda_M^n }}
\biggr),
\]
one can view this as the shrinking neighborhood around 0 that is being
treated as the region of inactive coefficients. The variance $\tau
_{1n}^{2}$ increases to $\infty$, where the rate depends on $p_n$.
However, when $p_n \prec\sqrt{n}$, $\tau_{1n}^{2}$ can be of constant
order [if $\lambda_m^n(\nu)$ is bounded away from zero].

Now consider the prior probability that a coefficient is nonzero
(denoted by $q_n$). The following calculation gives insight into the
choice of $q_n$. Let $K_n$ be a sequence going to $\infty$, then
\[
P \Biggl( \sum_{i =1}^{p_n}
Z_i> K_n \Biggr) \approx1 - \Phi \biggl(
\frac{K_n - p_n q_n}{\sqrt{p_n q_n (1-q_n)}} \biggr) \longrightarrow 0,
\]
if $p_n q_n$ is bounded. Therefore, we typically choose $q_n$ such that
$q_n \sim p^{-1}_n$. This can be viewed as a priori penalization of the
models with large size in the sense that the prior probability on
models with diverging number of covariates goes to zero. To this
respect, if $K$ is an initial upper bound for the size of the model
$t$, by choosing $q_n = c/p_n$ such that $\Phi ((K - c)/\sqrt{c}
) \approx1 - \alpha$, our prior probability on the models with
sizes greater than $K$ will be $\alpha$.

We would like to note that the hierarchical model considered by \citet
{George93} is similar to our model \eqref{modeleq}, but their prior
parameters are fixed and, therefore, do not satisfy our conditions. In
Section~\ref{ortho-sec}, we give an example illustrating model
selection inconsistency under fixed prior parameters.

\subsection{Methodology for variable selection}

We use the posterior distribution of the latent variables $Z_i$ to
select the active covariates. Note that the sample space of~$Z$,
denoted by $M$, has $2^{p_n}$ points, each of which corresponds to a
model. For this reason, we call $M$ the model space. To find the model
with the highest posterior probability is computationally challenging
for large $p_n$. In this paper, we use a simpler alternative, that is,
we use the $p_n$ marginal posterior probabilities $P(Z_i = 1|Y, X)$,
and select the covariates with the corresponding probability more than
a fixed threshold $\underbar{p} \in(0, 1)$. A threshold probability of
$0.5$ is a natural choice for $\underbar{p}$. This corresponds to what
\citet{Barbieri04} call the median probability model. In the orthogonal
design case, \citet{Barbieri04} showed that the median probability model
is an optimal predictive model. The median probability model may not be
the same as the maximum a posteriori (MAP) model in general, but the
two models are the same with probability converging to one under strong
selection consistency.

On the other hand, \citet{Dey08} argued that the median probability
model tends to underfit in finite samples. We also consider an
alternative by first ranking the variables based on the marginal
posterior probabilities and then using BIC to choose among different
model sizes. This option avoids the need to specify a threshold. In
either case, it is computationally advantageous to use the marginal
posterior probabilities, because we need fewer Gibbs iterations to
estimate only $p_n$ of them. The proposed methods based on marginal
posteriors achieve model selection consistency because the results in
Section~\ref{results-sec} assure that (i) the posterior probability of
the true model converges to~1, and (ii) the marginal posterior based
variable selection selects the true model with probability going to~1.
We now motivate these results and the necessity of sample size
dependent priors in a simple but illustrative case with orthogonal designs.

\section{Orthogonal design}\label{ortho-sec}
In this section, we consider the case where the number of covariates
$p_n < n$, and assume that the design matrix $X$ is orthogonal, that
is, $X'X = n I$. We also assume $\sigma^2$ to be known. Though this may
not be a realistic set-up, this simple case provides motivation for the
necessity of sample size dependent prior parameters as well as an
insight into the mechanism of model selection using these priors. At
this moment, we do not impose any assumptions on the prior parameters.
\textit{All the probabilities used in the rest of the paper are
conditional on~$X$}. Under this simple set-up, the joint posterior of
$\beta$ and $Z$ can be written as
\begin{eqnarray*}
&&P\bigl(\beta, Z \mid\sigma^2, Y
\bigr)
\\
&&\qquad\propto\exp \biggl\{{-\frac{1}{2 \sigma^2} \|Y - X\beta\| _2^2}
\biggr\} \prod_{i=1}^{p_n}
\bigl((1-q_n) \pi_0(\beta_i)
\bigr)^{1 - Z_i} \bigl(q_n \pi_1(
\beta_i) \bigr)^{Z_i}
\\
&&\qquad\propto\exp \biggl\{-\frac{1}{2 \sigma^2} \bigl(\beta'
X'X \beta- 2\beta 'X'Y\bigr) \biggr\}
\prod_{i=1}^{p_n} \bigl((1-q_n)
\pi_0(\beta_i) \bigr)^{1
- Z_i}
\bigl(q_n\pi_1(\beta_i)
\bigr)^{Z_i}
\\
&&\qquad\propto\exp \Biggl\{-\frac{n}{2 \sigma^2} \sum_{i=1}^p
(\beta_i - \hat{\beta}_i)^2 \Biggr\} \prod
_{i=1}^{p_n} \bigl((1-q_n)\pi
_0(\beta _i) \bigr)^{1 - Z_i}
\bigl(q_n\pi_1(\beta_i)
\bigr)^{Z_i},
\end{eqnarray*}
where for $k =0, 1$, $\pi_k(x) = \boldsymbol\phi(x, 0, \sigma^2
\tau
^2_{k,n})$ is the probability density function (p.d.f.) of the normal
distribution with mean zero and variance $\sigma^2 {\tau^2_{k,n}}$
evaluated at $x$, and $\hat{\beta}_i$ is the OLS estimator of $\beta
_i$, that is, $\hat{\beta}_i = X_i'Y/n$.

The product form of the joint posterior of $(Z_i, \beta_i)$ implies
that $(Z_i, \beta_i)$ and $\{(Z_j, \beta_j), j \neq i \}$ are
independent given data. Hence, the marginal posterior of $Z_i$ is given by
\[
P\bigl(Z_i \mid
\sigma^2, Y\bigr) \propto\int\exp \biggl\{-\frac{n}{2 \sigma
^2} (b -
\hat{\beta}_i)^2 \biggr\} \bigl((1-q_n)
\pi_0(b) \bigr)^{1 - Z_i} \bigl(q_n
\pi_1(b) \bigr)^{Z_i} \,db.
\]
Therefore,
%
\begin{equation}
\label{post-orth} %
P\bigl(Z_i = 0 \mid
\sigma^2, Y\bigr) = \frac{(1-q_n) E_{\hat{\beta}_i} (\pi_0
(B))}{(1-q_n)E_{\hat{\beta}_i} (\pi_0 (B)) +q_n E_{\hat{\beta}_i}
(\pi
_1 (B))},
\end{equation}
where $E_{\hat{\beta}_i}$ is the expectation under $B$ following the
normal distribution with mean $\hat{\beta}_i$ and variance $ \sigma^2/n
$. These expectations can be calculated explicitly, that is, for $k$ =
0 and 1,
\begin{eqnarray*}
\label{exp-orth} %
E_{\hat{\beta}_i} \bigl(
\pi_{k} (B) \bigr)& =& \frac{\sqrt {n}}{2\pi
\sigma\tau_{k,n}} \int\exp \biggl\{{-
\frac{n}{2 \sigma^2} (b- \hat {\beta}_i)^2} -
\frac{{b}^2}{2 {\tau^2_{k,n}}} \biggr\}\,db
\\
& =& \frac{1}{\sqrt{2 \pi}a_{k,n}} \exp \biggl\{-\frac{\hat{\beta
}^2_i}{2 a_{k,n}^2} \biggr\},
\end{eqnarray*}
where $a_{k,n} = \sqrt{\sigma^2/n + {\tau^2_{k,n}}}$.

This simple calculation gives much insight into the role of our priors
and the influence of the prior parameters on variable selection, which
we explain in some detail below. In the following subsections, we
assume that the $i$th covariate is identified as active if and only
if $P(Z_i = 1 \mid\sigma^2, Y) > 0.5$ for simplicity, and similar
arguments can be produced for threshold values other than 0.5.
\subsection{Fixed parameters}
Let us first consider the case of fixed parameters $\tau_{0n}^2 = \tau
_0^2 < \tau_{1n}^2 = \tau_1^2 $ and $q_n = q = 0.5$. We then have for
$k = 0, 1$,
%
\begin{equation}
E_{\hat{\beta}_i} \bigl(\pi_{k} (B)
\bigr) \cprob\frac{1}{{\tau_k}} \exp \biggl\{-\frac{\beta^2_i}{ 2 \tau^2_k}
\biggr\} \qquad\mbox{as } n \rightarrow\infty\mbox{ for } \beta_i \neq0.
\end{equation}
Now for $\beta_i = \tau_0 \neq0$, we have $ \exp \{- \beta
^2_i /
2 \tau^2_0  \}/\tau_0 > \exp \{-\beta^2_i/2 \tau^2_1
\}/\tau_1$ for any $\tau_1 \neq\tau_0$. Therefore, the limiting value
of $P(Z_i = 1 \mid\sigma^2, Y)$ will be less than 0.5 (with high
probability) as $n \rightarrow\infty$. This implies that even as $n
\rightarrow\infty$, we would not be able to identify the active
coefficient in this case.

\subsection{Shrinking \texorpdfstring{$\tau_{0,n}^2$}{tau 0,n 2}, fixed \texorpdfstring{$\tau_{1,n}^2$}{tau 1,n 2} and $q_n$}

Now consider the prior parameters such that $\tau_{1,n}^2$ and $q_n$ are
fixed, but $\tau_{0,n}^2$ goes to 0 with $n$. If $\beta_i = 0$,
$\sqrt {n} \hat{\beta}_i$ converges in distribution to the standard normal
distribution, and we have, for $k = 0, 1$,
\[
 \exp \biggl\{-\frac{\hat{\beta}^2_i}{2 (\sigma^2/n) + 2 \tau^2_{k,n}} \biggr\}
 =
O_P(1).
\]
In this case, \eqref{exp-orth} will imply that $E_{\hat{\beta}_i}
(\pi
_{1} (B)) = O_p(1)$, while $E_{\hat{\beta}_i} (\pi_{0} (B)) \cprob
\infty$. Therefore, from \eqref{post-orth}, we have $P(Z_i = 0 \mid
\sigma^2, Y) \cprob1$. For $\beta_i \neq0$, using $\hat{\beta}^2_i
\cprob{\beta}^2_i$ and the fact that $x e^{-rx^2}\rightarrow0$ as $x
\rightarrow\infty$ (for fixed $r >0$), we obtain $E_{\hat{\beta}_i}
(\pi_0 (B)) \rightarrow0$. As $E_{\hat{\beta}_i} (\pi_1 (B)) \sim c'$,
for some $c'>0$, we have $P(Z_i = 1 \mid\sigma^2, Y) \cprob1$.

To summarize, we have argued that $P(Z_i = 0 \mid\sigma^2, Y) \cprob
I{(\beta_i =0)}$, where $I(\cdot)$ is the indicator function.
That is, for orthogonal design matrices, the marginal posterior
probability of including an active covariate or excluding an inactive
covariate converges to one under shrinking prior parameter $\tau
_{0,n}^2$, with fixed parameters $\tau_{1,n}^2$ and $q_n$. However, it
should be noted that this statement is restricted to the convergence of
marginals of $Z$, and does not assure consistency of overall model
selection. To achieve this, we will need to allow $\tau_{1,n}^2$, $q_n$
to depend on the sample size, too.

\subsection{Shrinking and diffusing priors}
Note that the $i$th covariate is identified as active if and only if
\begin{eqnarray*}
&&P\bigl(Z_i = 1 \mid
\sigma^2, Y\bigr) > 0.5
\\
&&\qquad\Leftrightarrow q_n E_{\hat{\beta}_i} \bigl(\pi_1 (B)
\bigr) > (1-q_n) E_{\hat
{\beta}_i} \bigl(\pi_0 (B)\bigr)
\\
&&\qquad\Leftrightarrow\hat{\beta}^2_i \bigl(
a_{0,n}^{-2} - a_{1,n}^{-2} \bigr)> 2
\bigl(\log(1- q_n) a_{1,n} - \log q_n
a_{0,n} \bigr)
\\
&&\qquad\Leftrightarrow\hat{\beta}^2_i > 2 \bigl(\log(1-
q_n) a_{1,n} - \log q_n a_{0,n} \bigr)/
\bigl(a_{0,n}^{-2} - a_{1,n}^{-2}\bigr):=
\varphi_n.
\end{eqnarray*}

In particular, when $\tau^2_{0,n} = o(1/n)$, but the other parameters
$\tau^2_{1,n}$ and $q_n$ are fixed, we have $\varphi_n \sim\sigma^2
\log n /n $. Without loss of generality, assume that the first $|t|$
coefficients of $\beta$ are nonzero. For $i > |t|$, $\beta_i = 0$ which
implies that $n \hat{\beta}^2_i \eqdist\chi^2_1$. Therefore,
\begin{eqnarray*} P\biggl[\hat{\beta}^2_i
> \frac{\sigma^2 \log n }{n}\biggr] & = & P\bigl[\chi^2_1> \log
n\bigr]
\\
& \geq& \biggl(\frac{1}{\sqrt{\log n}} - \frac{1}{\sqrt{\log n}^3}\biggr)
e^{-{\log n}/{2}}
\\
& \geq& n^{-1/2 - \varepsilon},
\end{eqnarray*}
for $\varepsilon>0$ and sufficiently large $n$. Therefore, we have
\begin{eqnarray*} P\bigl[Z = t | \sigma^2, Y\bigr] &
\leq& P \biggl[\hat{\beta}^2_i \leq\frac
{\sigma
^2 \log n }{n},
\forall i > |t| \biggr]
\\
&\leq& \bigl(1 - n^{-1/2 - \varepsilon}\bigr)^{p_n - |t|}
\\
&\rightarrow& 0\qquad \mbox{if } p_n > n^{1/2 + 2\varepsilon}.
\end{eqnarray*}
The above argument shows that having $\tau^2_{1,n}$ and $q_n$ fixed
leads to inconsistency of selection if the number of covariates is much
greater than $\sqrt{n}$. In this case, the threshold $\varphi_n$ should
be larger to bound the magnitude of all the inactive covariates
simultaneously. By using the diffusing prior parameters Section~\ref{prior-sec}, the threshold will be $(2 + \delta) \sigma^2 \log p_n /n$
in place of $ \sigma^2 \log n/n $. Model selection consistency with
this threshold can be proved using similar arguments in the orthogonal
design case. We will defer the rigorous arguments to the next section.

\section{Main results}\label{results-sec}
In this section, we consider our model given by \eqref{modeleq} and
general design matrices. Because the model selection consistency holds
easily with $p_n = O(1)$, we assume throughout the paper that $p_n
\rightarrow\infty$ as $n \rightarrow\infty$.

\subsection{Conditions}
We first state the main conditions we use.
%
\begin{condition}[(On dimension $p_n$)] \label{dim-cond} $p_n = e^{n
d_n}$ for some $d_n\rightarrow0$ as $n \rightarrow\infty$, that is,
$\log p_n = o(n)$.
\end{condition}
%
\begin{condition}[(Prior parameters)]\label{prior-cond}
$n \tau_{0n}^2 = o(1)$, $ n \tau_{1n}^{2} \sim (n\vee p_n^{2 + 3
\delta}  )$, for some $\delta>0$, and $q_n \sim p^{-1}_n$.
\end{condition}
%
\begin{condition}[(On true model)]\label{model-cond}
$Y|X \sim N (X_t \beta_t + X_{ t^c} \beta_{t^c}, \sigma^2 I)$ where the
size of the true model $|t|$ is fixed. The coefficients corresponding
to the inactive covariates can be nonzero but satisfy $b_0:= \|X_{ t^c}
\beta_{t^c}\|_2 = O(1)$.
\end{condition}
For any fixed $K$, define
\[
\Delta_n(K):= \inf_{\{k: |k| < K |t|, k \not\supset t\}} \bigl\|(I - P_{k})
X_t \beta_t\bigr\|_2^2,
\]
where $P_{k}$ is the projection matrix onto the column space of $X_k$.
%
\begin{condition}[(Identifiability)] \label{identify-cond}
There is $K > 1 + 8/\delta$ such that $\Delta_n(K) > \gamma_n:= 5
\sigma
^2 |t| (1 + \delta) \log ( \sqrt{n} \vee p_n  )$. 
\end{condition}
%
%
\begin{condition}[(Regularity of the design)]\label{eigen-cond}
For some $\nu< \delta$, $\kappa< (K-1) \delta/2$,
\[
 \lambda_M^n \prec
\bigl( \bigl(n \tau_{0n}^2\bigr)^{-1} \wedge n
\tau_{1n}^2 \bigr)\quad \mbox{and}\quad \lambda_m^n
(\nu) \succeq \biggl( \frac{ n\vee p_n^{2
+ 2
\delta} }{n \tau_{1n}^{2}} \vee p_n^{- \kappa}
\biggr).
\]
\end{condition}
The moderateness of these conditions will be examined in some detail in
Section~\ref{cond-sec}.

\subsection{Results for fixed \texorpdfstring{$\sigma^2$}{sigma2}} We suppress $\nu$ and $K$
from the notation of $\lambda_m^n(\nu)$, $m_n(\nu)$ and $\Delta_n(K)$
for stating the results for convenience. In addition, we introduce the
following notation.
The posterior ratio of model $k$ with respect to the true model $t$ is
defined as
\[
\operatorname{PR}(k, t):={P\bigl(Z = k\mid Y,
\sigma^2\bigr)}/{P\bigl(Z = t\mid Y, \sigma^2\bigr)}.
\]
The following lemma gives an upper bound on the posterior ratio.
%
\begin{lem}\label{bf-lem}
Under Conditions \ref{prior-cond} and \ref{eigen-cond}, for any model
$k \neq t$, we have
\begin{eqnarray*}
\operatorname{PR}(k, t) &=& \frac{Q_k}{Q_t}
s_n^{|k| -|t|} \exp \biggl\{ -\frac{1}{2
\sigma^2} (
\tilde{R}_k - \tilde{R}_t) \biggr\}
\\
&\leq& w' \bigl(n \tau_{1n}^2
\lambda_m^n (1 - \phi_n)
\bigr)^ {-({1}/{2}) (r_{k}^* - r_t)} \bigl(\lambda_m^n
\bigr)^{- ({1}/{2}) |t \wedge
k^c|} s_n^{|k| -|t|}\\
&&{}\times\exp \biggl\{ -\frac{1}{2 \sigma^2} (\tilde{R}_k - \tilde
{R}_t) \biggr\},
\end{eqnarray*}
where $Q_k = | I + X D_k^{-1}X'|^{-1/2}$, $s_n = q_n/(1 - q_n) \sim
p_n^{-1}$, $w' > 0$ is a constant, $r_k = \operatorname{rank}( X_k)$, $r_k^* = r_k
\wedge m_n$, $\phi_n = o(1)$, $ \tilde{R}_k = Y' (I - X(D_k +
X'X)^{-1}X') Y $, and $D_k = \operatorname{Diag} (k \tau_{1n}^{-2} + (\mathbf{1} -k)
\tau_{0n}^{-2})$.
\end{lem}

The following arguments give some heuristics for the convergence of
pair-wise posterior ratio. Note that $\tilde{R}_k $ is the residual sum
of squares from a shrinkage estimator of $\beta$, and the term $LR_n:=
\exp\{ - (\tilde{R}_k - \tilde{R}_t)/2 \sigma^2 \}$ corresponds to the
usual likelihood ratio of the two models $k$ and $t$. Consider a model
$k$ that does not include one or more active covariates, then $ (\tilde
{R}_k - \tilde{R}_t)$ goes to $\infty$ at the same rate as $n$, because
it is (approximately) the difference in the residual sums of squares of
model $k$ and model $t$. We then have the posterior ratio converging to
zero since $LR_n \sim e^{-c n}$ for some $c>0$, and due to
Conditions \ref{dim-cond}--\ref{eigen-cond}, $P_n:= (n \tau_{1n}^2
\lambda_m^n (1 - \phi_n)  )^ {(r_t - r^*_k)/2} ({\lambda
_m^n})^{-|t \wedge k^c|/2} s_n^{|k| -|t|} (1 - \phi_n)^{-|t|/2} =
o(e^{cn})$. On the other hand, if the model $k$ includes all the active
covariates and one or more inactive covariates, we have $|k| > |t|$,
but $(\tilde{R}_k - \tilde{R}_t)$ is probabilistically bounded. The
posterior ratio in this case also converges to zero because $P_n$ goes to
zero. Note that when $r_k > r_t $, larger values of $\tau_{1n}^2$ will
imply smaller $P_n$. That is, the posterior ratio for large sized models
go to zero faster for larger values of $\tau_{1n}^2$. A~similar
observation is made by \citet{Ishwaran11}. To state our main result, we
first consider the posterior distributions of the models $Z$, assuming
the variance parameter $\sigma^2$ to be known. We consider the case
with the prior on $\sigma^2$ in Theorem~\ref{sigprior-thmm}.

\begin{thmm}\label{sigfix-thmm}
Assume Conditions \ref{dim-cond}--\ref{eigen-cond}. Under
model \eqref
{modeleq}, we have $P(Z = t \mid Y, \sigma^2) \cprob1$ as $n
\rightarrow\infty$, that is, the posterior probability of the true
model goes to 1 as the sample size increases to $\infty$.
\end{thmm}

\begin{remark} \label{sum-bf}
The statement of Theorem~\ref{sigfix-thmm} is equivalent to
%
\begin{equation}
\label{sum-bfeq} %
\frac{1 - P(Z = t\mid Y, \sigma^2)}{P(Z = t\mid Y, \sigma^2)} = \sum
_{k \neq t} \operatorname{PR}(k, t) \cprob0.
\end{equation}
\end{remark}

\begin{remark} It is worth noting that for Theorem~\ref{sigfix-thmm} to
hold, we do not actually need the true $\sigma^2$ to be known. Even for
a misspecified $\tilde{\sigma}^2 \neq\sigma^2$, $P(Z = t \mid Y,
\tilde
{\sigma}^2) \cprob1$ under the conditions $\Delta_n > \tilde{\sigma}^2
\gamma_n/ \sigma^2$ and $2 (1 +\delta) \tilde{\sigma}^2 > (2 +
\delta)
\sigma^2$. The same proof for Theorem~\ref{sigfix-thmm} works.
\end{remark}
To see why \eqref{sum-bfeq} holds, we provide specific rates of
convergence of individual posterior ratio summed over subsets of the
model space. We divide the set of models (excluding the model $t$) into
the following subsets:
\begin{longlist}[1.]
\item[1.]Unrealistically large models: $M_1 = \{k\dvtx r_k > m_n\}$, all the
models with dimension (i.e., the rank) greater than $ m_n$.
\item[2.]Over-fitted models: $M_2 = \{k\dvtx k \supset t, r_k \leq m_n\}$,
that is, the models of dimension smaller than $ m_n $ which include all
the active covariates plus one or more inactive covariates.
\item[3.]Large models: $M_3 = \{k\dvtx k \not\supset t, K |t| < r_k \leq
m_n \}$, the models which do not include one or more active covariates,
and dimension greater than $K |t|$ but smaller than $ m_n$.
\item[4.]Under-fitted models: $M_4 = \{k\dvtx k \not\supset t, r_k \leq K
|t| \}$, the models of moderate dimension which miss an active covariate.
\end{longlist}
The proof of Theorem~\ref{sigfix-thmm} shows the following results.
%
\begin{lem}[(Rates of convergence)]\label{rates-lem}
For some constants $c', w' >0 $ (which may depend on $\delta$), we have
\begin{longlist}[1.]
\item[1.] The sum of posterior ratio $\sum_{ k \in M_1}\operatorname{PR}(k, t) \preceq\exp
\{- w' n \}$, with probability at least $1 - 2 \exp\{ - c' n \}$.
\item[2.] The sum $\sum_{k \in M_2} \operatorname{PR}(k, t) \preceq v_n:=
(p_n^{-\delta/2}\wedge\frac{p_n^{1+\delta/2}}{ \sqrt{n}}  )$,
with probability greater than $1 - \exp\{ - c' \log p_n\}$.
\item[3.] The sum $\sum_{k \in M_3} \operatorname{PR}(k, t) \preceq\nu_n^{ (K-1) |t|/2
+1}$, with probability greater than $1 - \exp\{ - c' K |t| \log p_n\}$.
\item[4.] For some $w'' <1$, we have $\sum_{k \in M_4} \operatorname{PR}(k, t) \preceq
\exp\{ - w' (\Delta_n - w''\gamma_n)\}$, with probability greater than
$1 - \exp\{ - c' \Delta_n \}$.
\end{longlist}
\end{lem}

\subsection{Results with prior on \texorpdfstring{$\sigma^2$}{sigma2}}
We now consider the case with the inverse Gamma prior on the variance
parameter $\sigma^2$. Define the constant $w$ as $w:= \delta/8 (
1+\delta)^2$ in the rest of the section.
%
\begin{thmm}\label{sigprior-thmm}
Under the same conditions as in Theorem~\ref{sigfix-thmm}, if we only
consider models of dimension at most $ |t| + w n/\log p_n$, we have
$P(Z = t \mid Y) \cprob1$ as $n \rightarrow\infty$.
\end{thmm}

\begin{remark}
Note that the dimension of the models that need to be excluded for
Theorem~\ref{sigprior-thmm} to hold is in the order of $n/\log p_n$.
These are unrealistically large models that are uninteresting to us.
From now on, we implicitly assume this restriction when a prior
distribution is used for $\sigma^2$.
\end{remark}

The following corollary ensures that the variable selection procedure
based on the marginal posterior probabilities finds the right model
with probability tending to 1. It is a direct consequence of
Theorems \ref{sigfix-thmm} and \ref{sigprior-thmm}, but is particularly
useful for computations because it ensures that the marginal posterior
probabilities can be used for selecting the active covariates.
%
\begin{cor}\label{marg-cor}
Under the conditions of Theorem~\ref{sigprior-thmm}, we have for any $0
< \underbar{\textup{p}} < 1$,
$P [ P(Z_i = t_i \mid Y) > \underbar{\textup{p}} \mbox{ for all } i = 1,
\ldots, p_n  ] \rightarrow1$ as $n \rightarrow\infty$.
\end{cor}
\begin{pf}
Let $E_i$ be the event that the marginal posterior probability of
$i$th covariate $P(Z_i =t_i \mid Y) > \underbar{p}$. We shall show
that $P[ \bigcup_{i = 1}^{p_n} E_i^c] \rightarrow0$ as $n \rightarrow
\infty$. For each $i = 1, \ldots, p_n$, we have
\begin{eqnarray*}
P(Z_i \neq t_i \mid Y)
& = & \sum_{k:k_i \neq t_i} P(Z = k \mid Y)
\\
& \leq& \sum_{k \neq t} P(Z = k \mid Y)
\\
& = & 1 - P(Z = t|Y).
\end{eqnarray*}

Then $P[ \bigcup_{i = 1}^{p_n} E_i^c] = P [ P(Z_i = t_i \mid Y) \leq
\underbar{p} \mbox{ for some } i = 1, \ldots, p_n  ] \leq P [P(Z
= t|Y) \leq\underbar{p} ] \rightarrow0$, due to Theorem~\ref
{sigprior-thmm}.
\end{pf}

\section{Connection with penalization methods}\label{penalty-sec}
Due to Lemma~\ref{bf-lem}, the maximum aposteriori (MAP) estimate of
the model using our Bayesian set-up is equivalent to minimizing the
objective function
%
\begin{eqnarray}
\label{obj-fn} %
B(k) &:=& \tilde{R}_k
+ 2 \sigma^2 \bigl( - \bigl(|k| - |t|\bigr) \log s_n - \log
(Q_k /Q_t) \bigr)
\nonumber
\\[-8pt]
\\[-8pt]
\nonumber
& = & \tilde{R}_k + \bigl(|k| - |t|\bigr) \psi_{n,k},
\end{eqnarray}
where
\[
\psi_{n,k} = 2 \sigma^2
\biggl( - \log s_n - \frac
{\log ( Q_k /Q_t  )}{(|k| - |t|) } \biggr).
\]
Lemma~\ref{rates-lem} implies that with exponentially small
probability, the sum of posterior ratio of the models with dimension
greater than $ m_n $ goes to zero (exponentially) for the fixed $\sigma
$ case. We therefore focus on all the models with dimension less than
$m_n$ in this section. In addition, assume that the maximum and minimum
nonzero eigenvalues of models of size $2|t|$ are bounded away from
$\infty$ and 0, respectively. Then, due to Condition \ref{eigen-cond}
and the proof of Lemma~\ref{scale-lem}(iii), we have
%
\begin{equation}
\label{penal-bound} %
c \log(n \vee p_n)
\leq- \frac{\log ( Q_k/Q_t  )} {(r_k -
r_t) } \leq C \log(n \vee p_n),
\end{equation}
for some $0 < c \leq C < \infty$.

In particular, if the models with dimension less than $m_n$ are of full
rank, that is, $|k| = r_k$, then due to \eqref{penal-bound}, we have
%
\begin{equation}
\label{psi-bound} %
2 \sigma^2
c' \log(n \vee p_n) \leq\psi_{n,k} \leq2
\sigma^2 C' \log (n \vee p_n),
\end{equation}
where $0 < c' \leq C' < \infty$. As $ n \tau_{0n}^2 \lambda_{M}^n
\rightarrow0$, and $n \tau_{1n}^2 \lambda_{m}^n \rightarrow\infty$,
\[
\tilde{R}_k \sim Y' \bigl(I - X \bigl(1/
\tau_{1n}^2 + X' X\bigr)^{-1}X'
\bigr) Y = \| Y - \hat{Y}_{k}\|^2 + O(1).
\]
Therefore, the MAP estimate can be (asymptotically) described as the
model corresponding to minimizing the following objective function:
%
\begin{equation}
\label{l0} %
m (\beta):= \| Y - X \beta
\|_2^2 + \psi_{n,k} \bigl(\|\beta\|_0 -
|t|\bigr).
\end{equation}
Due to the bounds \eqref{psi-bound} on $\psi_{n,k}$, any inactive
covariate will be penalized in the order of $\log(n \vee p_n)$
irrespective of the size of the coeffecient. This is however not the
case with the $L_1$ penalty or SCAD penalty, which are directly
proportional to the magnitude of the coefficient in some interval
around zero.

The commonly used model selection criteria AIC and BIC are special
cases of $L_0$ penalization. The objective functions of AIC and BIC are
similar to $m(\beta)$, which have the quotient of penalty equal to $2$
and $\log n$ in place of $\psi_{n,k}$. Due to the results in Section~\ref{results-sec} and the above arguments, selection properties of our
proposed method are similar to those of the $L_0$ penalty. In
particular, it attempts to find the model with the least possible size
that could explain the conditional mean of the response variable.
A~salient feature of our approach is that the $L_0$-type penalization is
implied by the hierarchical model. The tuning parameters are more
transparent than those in penalization methods. Another feature to note
is that our model allows high (or even perfect) correlations among
inactive covatiates. This is practically very useful in high
dimensional problems because the number of inactive covariates is often
large and the singularity of the design matrix is a common occurrence.
Also, high correlations between active and inactive covariates is not
as harmful to the proposed method as they are to the $L_1$-type
penalties. This point is illustrated in Table~\ref{tab-highcor} of our
simulation studies in Section~\ref{simul-sec}.

\section{Discussion of the conditions}\label{cond-sec}
The purpose of this section is to demonstrate that Conditions~\ref
{dim-cond}--\ref{eigen-cond} that we use in Section~\ref{results-sec}
are quite mild. Condition~\ref{dim-cond} restricts the number of
covariates to be no greater than exponential in $n$, and Condition \ref
{prior-cond} provides the shrinking and diffusing rates for the spike
and slab priors, respectively. We note that Conditions \ref
{model-cond}--\ref{eigen-cond} allow $\beta$ to depend on $n$. For
instance, consider $p_n <n$ and the design matrix $X$ with $X'X/n
\rightarrow D$, where $D$ is a positive definite matrix. \citet
{Ishwaran05}, \citet{Zou06}, \citet{Bondell12} and \citet{Johnson12}
assumed this condition on the design under which Conditions \ref
{model-cond} and \ref{identify-cond} only require $\beta$ to be such that
\[
 \|\beta_{t^c}\|_2^2
= O \biggl(\frac{1}{n} \biggr) \quad\mbox{and}\quad \| \beta _t
\|_2^2 > c' \frac{\log n}{n},
\]
for some $c' >0$. Condition \ref{eigen-cond} is also satisfied in this
case, so Conditions \ref{model-cond}--\ref{eigen-cond} allow a wider
class of design matrices.

In general, Condition \ref{identify-cond} is a mild regularity
condition that allows us to identify the true model. It serves to
restrict the magnitude of the correlation between active and inactive
covariates, and also to bound the signal to noise ratio from below. The
following two remarks provide some insight into the role of Condition
\ref{identify-cond} in these aspects.
%
\begin{remark}\label{identify-rem}
Consider the case where the active coefficients $\beta_t$ are fixed.
We then have some $w' >0$, such that
\begin{eqnarray*}
\Delta_n (K) &\geq& \|
\beta_t\|_2^2 \inf_{\{k: |k| < K |t|, k \not
\supset t\}}
\phi_{\mathrm{min}} \bigl( X_t' (I - P_{k})
X_t \bigr)
\\
& \geq& w'n \inf_{\{k: |k| < K |t|, k \not\supset t\}} \phi_{\mathrm{min}}
\biggl(\frac{ X_{k \vee t}' X_{k \vee t}}{n} \biggr),
\end{eqnarray*}
where we have used the fact that $\phi_{\mathrm{min}}  ( X_{k \vee t}' X_{k
\vee t} ) \leq\phi_{\mathrm{min}}  ( X_t' (I - P_{k}) X_t
)$. To
see this, we just need to consider the cases where $X_{k \vee t}$ is of
full rank. Then it follows from the observation that $ ( X_t' (I -
P_{k}) X_t  )^{-1}$ is a submatrix of $ ( X_{k \vee t}' X_{k
\vee t} )^{-1}$.
Therefore, Condition \ref{identify-cond} is satisfied if the minimum
eigenvalues of the submatrices of $X'X/n$ with size smaller than $(K+1)
|t|$ are uniformly larger than $c' \log(n \vee p_n)/n$. In the other
end of the spectrum, where the inactive covariates can be perfectly
correlated, Condition \ref{identify-cond} could still hold.
\end{remark}

\begin{remark} If the infimum of $\phi_{\mathrm{min}}  ( X_t' (I - P_{k})
X_t /n )$ is uniformly bounded away from zero, then $\Delta_n (K)
\geq w' n \|\beta_t\|_2^2$. Then Condition \ref{identify-cond} is
satisfied if
\[
\biggl\llVert \frac{\beta_t}{\sigma}\biggr\rrVert _2^2 \geq
\frac{c' \log(n
\vee p_n)}{n}.
\]
\end{remark}

Condition \ref{eigen-cond} provides conditions on the eigenvalues of
the Gram matrix in terms of the prior parameters. The condition is
weaker than the assumption that the maximum and minimum nonzero
eigenvalues of the Gram matrix are bounded away from infinity and zero,
respectively. In Condition \ref{eigen-cond}, $\lambda_M^n \prec(n
\tau
_{0n}^2)^{-1}$ will be satisfied if $\tau_{0n}^2$ is small enough.
However, the assumption on $\lambda_m^n (\nu) $ is nontrivial as it
needs to be greater than $p_n^{-\kappa}$. We now show that this
requirement is satisfied with high probability if the design matrix
consists of independent sub-Gaussian rows.

\begin{lem}[(MNEV for sub-Gaussian random matrices)] \label{mnev-lemma}
Suppose that the rows of $X_{n \times p_n}$ are independent isotropic
sub-Gaussian random vectors in $R^{p_n}$. Then there exists a $\nu>0$
such that, with probability greater than $1 - \exp(- w' n)$,
\[
\inf_{|k| \leq m_n(\nu)}
\phi_{\mathrm{min}} \biggl(\frac{X_k' X_k}{n} \biggr) >0.
\]
\end{lem}

A proof of Lemma~\ref{mnev-lemma} is provided in Section~\ref{proofs-sec}. Lemma~\ref{mnev-lemma} implies that the Gram matrix of a
sub-Gaussian design matrix has the minimum eigenvalues of all the
$m_n(\nu)$ dimensional\vadjust{\goodbreak} submatrices to be uniformly bounded away
from
zero. This clearly is stronger than Condition \ref{eigen-cond}, which
only requires the minimum nonzero eigenvalues to be uniformly greater
than $p_n ^{-\kappa}$. In particular, unlike the restricted isometry
conditions which control the minimum eigenvalue, Condition \ref
{eigen-cond} allows the minimum eigenvalue to be exactly zero to allow
even perfect correlation among inactive (or active) covariates.

\section{Computation} \label{compu-sec}

The implementation of our proposed method involves using the Gibbs
sampler to draw samples from the posterior of $Z$. The full
conditionals are standard distributions due to the use of conjugate
priors. The conditional distribution of $\beta$ is given by
\[
f\bigl(\beta\mid Z, \sigma^2, Y
\bigr) \propto\exp \biggl\{{-\frac{1}{2 \sigma^2
} \|Y - X\beta\|_2^2}
\biggr\} \prod_{i=1}^{p_n} \phi \bigl(
\beta_i, 0, \sigma^2 \tau^2_{Z_i,n}
\bigr),
\]
where $ \phi(x, 0, \tau^2)$ is the p.d.f. of the normal distribution with
mean zero, and variance $\tau^2$ evaluated at $x$.
This can be rewritten as
\[
f\bigl(\beta\mid Z = k, \sigma^2, Y
\bigr) \propto\exp \biggl\{-\frac{1}{2
\sigma
^2} \bigl(\beta'
X' X \beta- 2\beta'X'Y\bigr) \biggr\}
\exp \biggl\{-\frac{1}{2
\sigma^2} \beta' D_k \beta
\biggr\},
\]
where $D_k = \operatorname{Diag} (\tau_{k_i, n}^{-2})$. Hence, the conditional
distribution of $\beta$ is given by $\beta\sim N(m, \sigma^2 V)$,
where $V =(X' X + D_k)^{-1}$, and $ m = V X'Y$. Furthermore, the
conditional distribution of $Z_i$ is
%
\[
P\bigl(Z_i = 1\mid\beta, \sigma^2\bigr) =
\frac{q_n \phi(\beta_i, 0, \sigma^2
{\tau^2_{1,n}} )}{q_n \phi(\beta_i, 0, \sigma^2 {\tau^2_{1,n}} ) +
(1-q_n) \phi(\beta_i, 0, \sigma^2 {\tau^2_{0,n}} ) }.
\]
The conditional of $\sigma^2$ is the inverse Gamma distribution $\operatorname{IG}(a,
b)$ with $a = \alpha_1 + n/2 + p_n/2$, and $b = \alpha_2 + \beta'D_k
\beta/2 + (Y - X \beta)' (Y - X\beta)/2$.

The only possible computational difficulty in the Gibbs sampling
algorithm is the step of drawing from the conditional distribution of
$\beta$, which is a high dimensional normal distribution for large
values of $p_n$. However, due to the structure of the covariance matrix
$(X' X + D_k)^{-1}$, it can be efficiently sampled using block updating
that only requires drawing from smaller dimensional normal
distributions. Details of the block updating can be found in \citet{Ishwaran05}.

\section{Simulation study}\label{simul-sec}
In this section, we study performance of the proposed method in several
experimental settings, and compare it with some existing variable
selection methods. We will refer to the proposed method as BASAD for
BAyesian Shrinking And Diffusing priors.

The proposed BASAD method has three tuning parameters. In all our
empirical work, we use
\[
\tau_{0n}^2=
\frac{\hat\sigma^2 }{10n},    \qquad \tau_{1n}^2= \hat
\sigma^2 \max \biggl( \frac{p_n^{2.1}}{100 n}, \log n \biggr),
\]
where $\hat\sigma^2$ is the sample variance of $Y$, and we choose $q_n
= P[Z_i = 1]$ such that $P  [\sum_{i=1}^{p_n} Z_i = 1 >K  ] =
0.1$, for a prespecified value of $K$. Our default value is $K=\max(10
, \log(n))$, unless otherwise specified in anticipation of a less
sparse model. The purpose of using $\hat\sigma^2$ is to provide
appropriate scaling. If a preliminary model is available, it is better
to use as $\hat\sigma^2$ the residual variance from such a model. It is
clear that those choices are not optimized for any given problem, but
they provide a reasonable assessment on how well BASAD can do. In the
simulations, we use 1000 burn-in iterations for the Gibbs sampler
followed by 5000 updates for estimating the posterior probabilities. As
mentioned in Section~\ref{model-sec}, we consider both the median
probability model (denoted by BASAD) and the BIC-based model (denoted
by BASAD.BIC) where the threshold for marginal posterior probability is
chosen by the BIC.
The R function used for obtaining the results in this section is publicly
available on the authors' website.

In this paper, we report our simulation results for six cases under
several $(n, p)$ combinations, varied correlations, signal strengths
and sparsity levels.

\begin{itemize}
\item \textit{Case} 1: In the first case, we use the set-up of \citet{Johnson12}
with $p = n$. Two sample sizes, $n = 100$ and $n = 200$, are
considered, and the covariates are generated from the multivariate
normal distributions with zero mean and unit variance. The compound
symmetric covariance with pairwise covariance of $\rho= 0.25$ is used
to represent correlation between covariates. Five covariates are taken
active with coefficients $\beta_t = (0.6, 1.2, 1.8, 2.4, 3.0)$. This is
a simple setting with moderate correlation between covariates and
strong signal strength.

\item \textit{Case} 2: We consider the $p >n$ scenario with $(n, p) = (100,
500)$ and $(n, p) = (200, 1000)$, but the other parameters are same as
in case 1.

For the next three cases (cases 3--5), we keep $(n,p) = (100,500)$ but
vary model sparsity, signal strength and correlation among covariates.

\item \textit{Case} 3: We keep $\rho= 0.25 $ and $|t|=5$ but have low signals
$\beta_t = (0.6, 0.6, 0.6,\break   0.6, 0.6)$.

\item \textit{Case} 4: We consider a block covariance setting where the active
covariates have common correlation ($\rho_1$) equal to 0.25, the
inactive covariates have common correlation ($\rho_3$) equal to 0.75
and each pair of active and inactive covariate has correlation ($\rho
_2$) 0.50. The other aspects of the model are the same as in case 1.

\item \textit{Case} 5: We consider a less sparse true model with $|t| = 25$ and
$\beta_t$ is the vector containing 25 equally spaced values between 1
and 3 (inclusive of 1 and 3).

\item \textit{Case} 6: We consider the more classical case of $n > p$ with
$(n,p) = (100,50)$ and $(n,p) = (200,50)$. Following \citet{Bondell12},
the covariates are drawn from a normal distribution with the covariance
matrix distributed as the Wishart distribution centered at the identity
matrix with $p$ degrees of freedom. Three of the 50 covariates are
taken to be active with their coefficients drawn from the uniform
distribution $U(0, 3)$ to imply a mix of weak and strong signals.
\end{itemize}

\begin{table}
\caption{Performance of BASAD for case 1:
$n = p$. The methods under comparison are $\mathit{piMOM}$ with nonlocal priors
of
\citet{Johnson12}, $\mathit{BCR.Joint}$ of \citet{Bondell12}, $\mathit{SpikeSlab}$ of
\citet{Ishwaran05} and three penalization methods Lasso, elastic net
(EN), and SCAD tuned by the BIC. The other columns of the table are as
follows: $pp_0$ and $pp_1$ (when applicable) are the average posterior
probabilities of inactive and active variables, respectively; $Z=t$ is
the proportion that the exact models is selected; $Z \supset t$ is the
proportion that the selected model contains all the active covatiates;
FDR is the false discovery rate, and $\mathit{MSPE}$ is the mean squared
prediction error of the selected models}
\label{tab1} 
%
\begin{tabular*}{\textwidth}{@{\extracolsep{\fill}}lcccccc@{}}
\hline
& $\bolds{pp_0}$ & $\bolds{pp_1}$ & $ \bolds{Z = t} $& $ \bolds{Z \supset t} $
& \textbf{FDR}& \textbf{MSPE}\\ %
\hline
&\multicolumn{6}{c}{$(n, p) = (100, 100); \rho= 0.25; |t| = 5 $} \\
BASAD&0.016&0.985&0.866&0.954&0.015&1.092\\
BASAD.BIC&0.016&0.985&0.066&0.996&0.256&1.203\\
piMOM&0.012&0.991&0.836&0.982&0.030&1.083\\
BCR.Joint&&&0.442&0.940&0.157&1.165\\
SpikeSlab&&&0.005&0.216&0.502&1.660\\
Lasso.BIC&&&0.010&0.992&0.430&1.195\\
EN.BIC&&&0.398&0.982&0.154&1.134\\
SCAD.BIC&&&0.356&0.990&0.160&1.157\\[3pt]
&\multicolumn{6}{c}{$(n, p) = (200, 200); \rho= 0.25; |t| = 5 $} \\
BASAD&0.002&1.000&0.944&1.000&0.009&1.037\\
BASAD.BIC&0.002&1.000&0.090&1.000&0.187&1.087\\
piMOM&0.003&1.000&0.900&1.000&0.018&1.038\\
BCR.Joint&&&0.594&0.994&0.102&1.064\\
SpikeSlab&&&0.008&0.236&0.501&1.530\\
Lasso.BIC&&&0.014&1.000&0.422&1.101\\
EN.BIC&&&0.492&1.000&0.113&1.056\\
SCAD.BIC&&&0.844&1.000&0.029&1.040\\
\hline
\end{tabular*}
\end{table}

\begin{table}
\caption{Performance of BASAD for case 2: $p >
n$}\label{tab2} 
\begin{tabular*}{\textwidth}{@{\extracolsep{\fill}}lcccccc@{}}
\hline
& $\bolds{pp_0}$ & $\bolds{pp_1}$ & $ \bolds{Z = t} $& $ \bolds{Z \supset t} $
& \textbf{FDR}& \textbf{MSPE}\\ %
\hline
&\multicolumn{6}{c}{$(n, p) = (100, 500); \rho= 0.25; |t| = 5$} \\
BASAD&0.001&0.948&0.730&0.775&0.011&1.130\\
BASAD.BIC&0.001&0.948&0.190&0.915&0.146&1.168\\
BCR.Joint&&&0.070&0.305&0.268&1.592\\
SpikeSlab&&&0.000&0.040&0.626&3.351\\
Lasso.BIC&&&0.005&0.845&0.466&1.280\\
EN.BIC&&&0.135&0.835&0.283&1.223\\
SCAD.BIC&&&0.045&0.980&0.328&1.260\\[3pt]
%
&\multicolumn{6}{c}{$(n, p) = (200, 1000); \rho= 0.25; |t| = 5$} \\
BASAD&0.000&0.986&0.930&0.950&0.000&1.054\\
BASAD.BIC&0.000&0.986&0.720&0.990&0.046&1.060\\
BCR.Joint&&&0.090&0.250&0.176&1.324\\
SpikeSlab&&&0.000&0.050 &0.574 &1.933\\
Lasso.BIC&&&0.020&1.000&0.430&1.127\\
EN.BIC&&&0.325&1.000&0.177&1.077\\
SCAD.BIC&&&0.650&1.000&0.091&1.063\\
\hline
\end{tabular*}
\end{table}

\begin{table}[b]
\caption{Performance of BASAD for case 3: $(n, p) =
(100, 500)$}\label{tablowsignal} 
\begin{tabular*}{\textwidth}{@{\extracolsep{\fill}}lcccccc@{}}
\hline
& $\bolds{pp_0}$ & $\bolds{pp_1}$ & $ \bolds{Z = t} $& $ \bolds{Z \supset t} $
& \textbf{FDR}& \textbf{MSPE}\\ %
\hline
&\multicolumn{6}{c}{$ \rho= 0.25; |t| = 5; \beta_t = (0.6,
0.6,0.6,0.6,0.6)$} \\
BASAD&0.002&0.622&0.185&0.195&0.066&2.319\\
BASAD.BIC&0.002&0.622&0.160&0.375&0.193&1.521\\
BCR.Joint&&&0.030&0.315&0.447&1.501\\
SpikeSlab&&&0.000&0.000&0.857&2.466\\
Lasso.BIC&&&0.000&0.520&0.561&1.555\\
EN.BIC&&&0.040&0.345&0.478&1.552\\
SCAD.BIC&&&0.045&0.340&0.464&1.561\\
\hline
\end{tabular*}
\end{table}

\begin{table}
\caption{Performance of BASAD for case 4: $(n, p) =
(100, 500)$}\label{tab-highcor} 
\begin{tabular*}{\textwidth}{@{\extracolsep{\fill}}lcccccc@{}}
\hline
& $\bolds{pp_0}$ & $\bolds{pp_1}$ & $ \bolds{Z = t} $& $ \bolds{Z \supset t} $
& \textbf{FDR}& \textbf{MSPE}\\ %
\hline
&\multicolumn{6}{c}{$\rho_1 = 0.25, \rho_2 = 0.50, \rho_3 = 0.75$}
\\
BASAD&0.002&0.908&0.505&0.530&0.012&\phantom{0}1.199\\
BASAD.BIC&0.002&0.908&0.165&0.815&0.179&\phantom{0}1.210\\
BCR.Joint&&&0.000&0.000&0.515&\phantom{0}2.212\\
SpikeSlab&&&0.000&0.000&0.995&10.297\\
Lasso.BIC&&&0.000&0.015&0.869&\phantom{0}8.579\\
EN.BIC&&&0.000&0.000&0.898&\phantom{0}8.360\\
SCAD.BIC&&&0.000&0.000&0.899&\phantom{0}8.739\\
\hline
\end{tabular*}
\end{table}

\begin{table}[b]
\caption{Performance of BASAD for case 5: $(n, p) =
(100, 500)$. In this case, two versions of BASAD are included, where
BASAD.K10 uses our default value of $K=10$, and BASAD.K50 uses a less
sparse specification of $K=50$}\label{tabnz25} 
\begin{tabular*}{\textwidth}{@{\extracolsep{\fill}}lcccccd{3.3}@{}}
\hline
& $\bolds{pp_0}$ & $\bolds{pp_1}$ & $ \bolds{Z = t} $& $ \bolds{Z \supset t} $
& \textbf{FDR}& \multicolumn{1}{c@{}}{\textbf{MSPE}}\\ %
\hline
&\multicolumn{6}{c}{$ \rho= 0.25; |t| = 25 $} \\
BASAD.K50&0.020&0.988&0.650&0.950&0.036&3.397\\
BASAD.BIC.K50&0.020&0.988&0.005&0.960&0.283&4.019\\
BASAD.K10&0.003&0.548&0.405&0.420& 0.011&170.862\\
BASAD.BIC.K10& 0.003&0.548&0.035&0.430&0.076&88.881\\
BCR.Joint&&&0.000&0.000&0.622&49.299\\
SpikeSlab&&&0.000&0.000&0.816&111.911\\
Lasso.BIC&&&0.000&0.005&0.685&58.664\\
EN.BIC&&&0.000&0.000&0.693&59.058\\
SCAD.BIC&&&0.000&0.000&0.666&72.122\\[3pt]
%
&\multicolumn{6}{c}{$ \rho= 0.75; |t| = 25 $} \\ [0.5ex]
BASAD.K50&0.048&0.914&0.005&0.355&0.289&6.103\\
BASAD.BIC.K50&0.048&0.914&0.000&0.445&0.498&6.611\\
BASAD.K10&0.003&0.298&0.025&0.030&0.018&349.992\\
BASAD.BIC.K10&0.003&0.298&0.000&0.060&0.087&61.709\\
BCR.Joint&&&0.000&0.000&0.772&34.113\\
SpikeSlab&&&0.000&0.000&0.899&48.880\\
Lasso.BIC&&&0.000&0.000&0.734&24.310\\
EN.BIC&&&0.000&0.000&0.754&29.171\\
SCAD.BIC&&&0.000&0.000&0.736&27.236\\
\hline
\end{tabular*}
\end{table}

\begin{table}
\caption{Performance of BASAD for case 6: $n >
p$}\label{tabnmore}
\centering
\begin{tabular*}{\textwidth}{@{\extracolsep{\fill}}lcccccc@{}}
\hline
& $\bolds{pp_0}$ & $\bolds{pp_1}$ & $ \bolds{Z = t} $& $ \bolds{Z \supset t} $
& \textbf{FDR}& \textbf{MSPE}\\ %
\hline
&\multicolumn{6}{c}{$(n,p) = (100,50)$} \\
BASAD&0.037&0.899&0.654&0.714&0.026&1.086\\
BASAD.BIC&0.037&0.899&0.208&0.778&0.267&1.151\\
piMOM&0.011&0.892&0.656&0.708&0.021&1.066\\
SpikeSlab&&&0.064&0.846&0.567&1.226\\
BCR.Joint&&&0.336&0.650&0.216&1.124\\
Lasso.BIC&&&0.076&0.744&0.397&1.152\\
EN.BIC&&&0.378&0.742&0.194&1.110\\
SCAD.BIC&&&0.186&0.772&0.284&1.147\\[3pt]
&\multicolumn{6}{c}{$(n,p) = (200,50)$} \\
BASAD&0.026&0.926&0.738&0.784&0.017&1.029\\
BASAD.BIC&0.026&0.926&0.338&0.842&0.193&1.055\\
piMOM&0.005&0.908&0.694&0.740&0.020&1.036\\
BCR.Joint&&&0.484&0.770&0.133&1.045\\
SpikeSlab&&&0.038&0.900&0.629&1.121\\
Lasso.BIC&&&0.082&0.752&0.378&1.059\\
EN.BIC&&&0.428&0.748&0.165&1.039\\
SCAD.BIC&&&0.358&0.812&0.193&1.046\\
\hline
\end{tabular*}\vspace*{-5pt}
\end{table}

The summary of our results are presented in Tables \ref{tab1}--\ref
{tabnmore}. In those tables, BASAD denotes the median probability
model, BASAD.BIC denotes the model obtained by using the threshold
probability chosen by the BIC. Three competing Bayesian model selection
methods are: (1) piMOM, the nonlocal prior method proposed by \citet
{Johnson12} but only when $p \le n$; (2) BCR.Joint, the Bayesian joint
credible region method of \citet{Bondell12} (using the default priors
followed by an application of BIC); (3) SpikeSlab, the generalized
elastic net model obtained using the R package spikeslab [\citet
{Ishwaran10}] for the spike and slab method of \citet{Ishwaran05}. Three
penalization methods under consideration are: (1) LASSO; (2) Elastic
Net (EN); and (3) SCAD, all tuned by the BIC. Our simulation experiment
used 500 data sets from each model when $n \geq p$, but used 200 data
sets when $p > n$ to aggregate the results.

The columns of the tables show the average marginal posterior
probability assigned to inactive covariates and active covariates
($pp_0$ and $pp_1$, resp.), proportion of choosing the true
model ($ Z = t $), proportion of including the true model ($ Z \supset
t $) and false discovery rate (FDR). The last column (MSPE) gives the
average test mean squared prediction error based on $n$ new
observations as testing data. From our simulation experiment, we have
the following findings:
\begin{longlist}[(iii)]
\item[(i)] The Bayesian model selection methods BASAD and piMOM (whenever
available) tend to perform better then the other methods in terms of
selecting the true model and controlling the false discovery rate in
variable selection, and our proposed BASAD stands out in this regard.
The penalization methods often have higher probabilities of selecting
all the active covariates at the cost of overfitting and false
discoveries. In terms of the prediction error, however, BASAD does not
always outperform its competitors, but remains competitive.

\item[(ii)] When the signals are low (case 3), all the methods under
consideration have trouble finding the right model, and BASAD.BIC
results in lower prediction error than BASAD with 0.5 as the threshold
for posterior probabilities.
In most cases, BASAD.BIC leads to slightly higher false positive rates
than BASAD with similar prediction errors.

\item[(iii)] In case 4, there is a moderate level of correlation among
inactive covariates and some level of correlation between active and
inactive covariates. This is where BASAD outperforms the other methods
under consideration because BASAD is similar to the $L_0$ penalty and
is able to accommodate such correlations well. Please refer to our
discussion in Sections \ref{penalty-sec} and \ref{cond-sec}.

\item[(iv)] When the true model is not so sparse and has $|t|= 25$ active
covariates (case 5), our default choice of $K=10$ in BASAD did not
perform well, which is not surprising. In fact, no other methods under
consideration did well in this case, highlighting the difficulty of
finding a nonsparse model with a limited sample size. On the other
hand, there is some promising news. If we anticipate a less sparse
model with $K=50$, the proposed method BASAD improved the performance
considerably. Our empirical experience suggests that if we are
uncertain about the level of sparsity of our model, we may use a
generous choice of $K$ or use BIC to choose between different values of $K$.\vadjust{\goodbreak}
\end{longlist}

\section{Real data example}\label{data-sec}

In this section, we apply our variable selection method to a real data
set to examine how it works in practice. We consider the data from an
experiment conducted by \citet{Lan06} to study the genetics of two
inbred mouse populations (B6 and BTBR). The data include expression
levels of 22,575 genes of 31 female and 29 male mice resulting in a
total of 60 arrays. Some physiological phenotypes, including the
numbers of phosphoenopyruvate carboxykinase (PEPCK) and
glycerol-3-phosphate acyltransferase (GPAT) were also measured by
quantitative real-time PCR. The gene expression data and the phenotypic
data are available at GEO (\url{http://www.ncbi.nlm.nih.gov/geo}; accession
number GSE3330). \citet{Zhang09} used orthogonal components regression
to predict each phenotype based on the gene expression data. \citet
{Bondell12} used the Bayesian credible region method for variable
selection on the same data.

Because this is an ultra-high dimensional problem with $p_n = 22\mbox{,}575$,
we prefer to perform simple screenings of the genes first based on the
magnitude of marginal correlations with the response. The power of
marginal screening has been recognized by \citet{Fan08}. After the
screening, the dataset for each of the responses consisted of $p = 200$
and $400$ predictors (including the intercept and gender) by taking 198
and 398 genes based on marginal screening. We performed variable
selection with BASAD along with LASSO, SCAD and the BCR method.
Following \citet{Bondell12}, we randomly split the sample into a
training set of 55 observations and a test set with the remaining five
observations. The fitted models using the training set were used to
predict the response in the test set. This process was repeated 100
times to estimate the prediction error.

\begin{figure}

\includegraphics{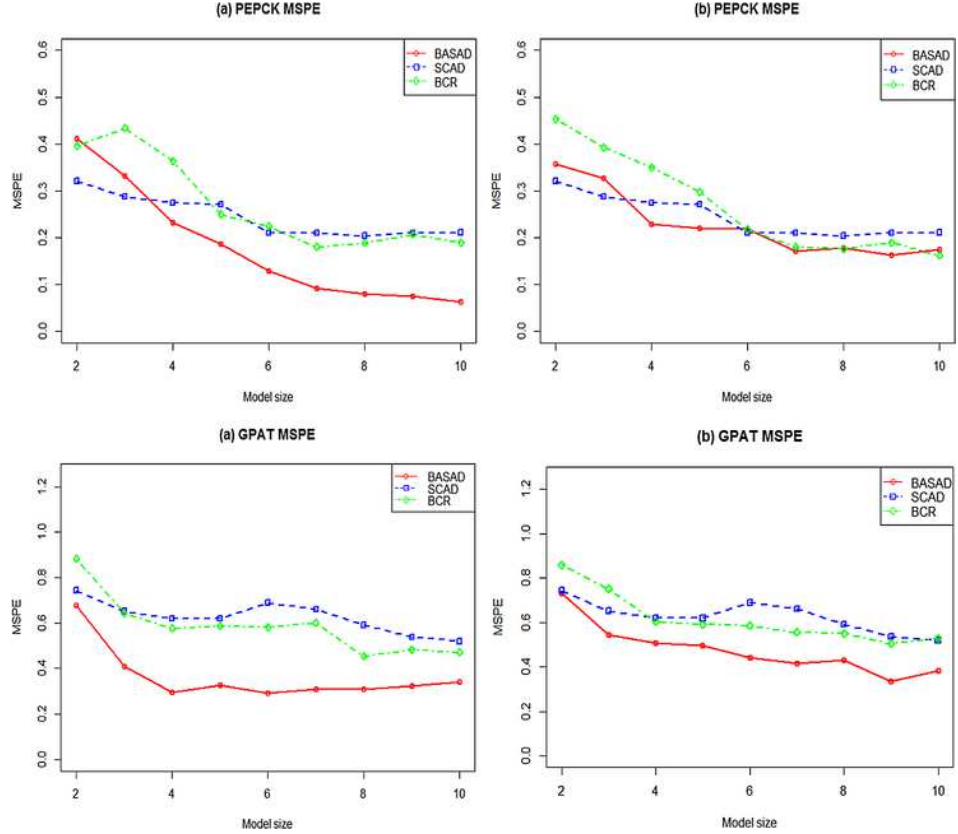}

\caption{Mean squared prediction error (MSPE) versus
model size for analyzing PEPCK and GPAT in the upper and lower panel,
respectively, \textup{(a)} $p = 200$ and \textup{(b)} $p = 400$.}
\label{fig1}
\end{figure}

In Figure~\ref{fig1}, we plot the average mean square prediction error
(MSPE) for models of various sizes chosen by BASAD, BCR and SCAD
methods for the two responses PEPCK and GPAT. We find that the MSPE of
BASAD is mostly smaller than that for other methods across different
model sizes. In particular, BASAD chooses less correlated variables and
achieves low MSPE with fewer predictive genes than the other methods.
We also note that the 10-covariate models chosen by BASAD is very
different (with the overlap of just one covariate for PEPCK and three
covariates for GPAT) from those of SCAD which chose mostly the same
covariates as LASSO. There are four common covariates identified by
both BASAD and BCR methods. When we perform a linear regression by
including the covariates chosen by BASAD and SCAD, we noticed that
majority of the covariates chosen by BASAD are significant, which
indicates that those genes chosen by BASAD are significant in
explaining the response even in the presence of those chosen using
SCAD. Most of the genes selected by SCAD, however, are not significant
in the presence of those chosen by BASAD. Despite the evidence in favor
of the genes selected by BASAD in this example, we must add that the
ultimate assessment of a chosen model would need to be made by
additional information from the subject matter science and/or
additional experiment.

\section{Conclusion}\label{conclusion-sec}

In this paper, we consider a Bayesian variable selection method for
high dimensional data based on the spike and slab priors with shrinking
and diffusing priors. We show under mild conditions that this approach
achieves strong selection consistency in the sense that the posterior
probability of the true model converges to one. The tuning parameters
needed for the prior specifications are transparent, and a standard
Gibbs sampler can be used for posterior sampling. We also provide the
asymptotic relationship between the proposed approach and the $L_0$
penalty for model selection. Simulation studies in Section~\ref{simul-sec} and real data example in Section~\ref{data-sec} show
evidence that the method performs well in a variety of settings even
though we do not attempt to optimize the tuning parameters in the
proposed method.

The strong selection consistency of Bayesian methods has not been
established in the cases of $p >n$ until very recently. For higher
dimensional cases, we just became aware of \citet{Liang13}, which
provided the strong selection consistency for Bayesian subset selection
based on the theory developed by \citet{Jiang07} for posterior density
consistency. However, to translate density consistency into selection
consistency, \citet{Liang13} imposed a condition on the posterior
distribution itself, which is not verifiable directly. The techniques
we use in this paper might also be used to complete the development of
their theory on strong selection consistency.

Throughout the paper, we assume Gaussian errors in the regression
model, but this assumption is not necessary to obtain selection
consistency. For proving Lemma~\ref{bf-lem}, we did not need
assumptions on the error distribution, and to prove Theorem~\ref
{sigprior-thmm}, we just need deviation inequalities of the quadratic
forms $\varepsilon' P_k \varepsilon$, which follow the chi-squared
distribution for normal errors. Similar proofs with an application of
deviation inequalities for other error distributions would work. For
instance, \citet{Hsu12} provide deviation inequalities for quadratic
forms of sub-Gaussian random variables.

The primary focus of our paper is model selection consistency. The
model is selected by averaging over the latent indicator variables
drawn from the posterior distributions. The strengths of different
model selection methods need to be evaluated differently if prediction
accuracy is the goal. In our empirical work, we have included
comparisons of the mean squared prediction errors, and found that our
proposed method based on default tuning parameters is highly
competitive in terms of prediction. However, improvements are possible,
mainly in the cases of low signals, if the parameters are tuned by BIC
or cross-validation, or if model-averaging is used instead of the
predictions from a single model.

\section{Proofs}\label{proofs-sec}
In this section, we prove Lemmas \ref{bf-lem} and  \ref{mnev-lemma}.
Please refer to \citet{Narisetty14} for proofs of the remaining results.
\begin{pf*}{Proof of Lemma~\ref{bf-lem}}
The joint posterior of $\beta, \sigma^2, Z$ under model \eqref{modeleq}
is given by
%
\begin{eqnarray}
\label{fullpost-eq}&& P\bigl(\beta, Z =k, \sigma^2 \mid Y\bigr)
\nonumber
\\
&&\qquad\propto\exp \biggl\{ -\frac{1}{2
\sigma^2} \bigl(\|Y - X\beta\|_2^2
- \beta' D_k \beta- 2\alpha_2 \bigr)
\biggr\} \sigma^{- 2 ({n}/{2} + {p_n}/{2} +
\alpha_1 +1 )}\\
&&\qquad\quad{}\times |D_k|^{{1}/{2}}
s_n^{|k|},\nonumber
\end{eqnarray}
where $D_k = \operatorname{Diag} (k \tau_{1n}^{-2} + (\mathbf{1} -k) \tau_{0n}^{-2})$,
$s_n =q_n/(1-q_n)$, $\alpha_1, \alpha_2$ are the parameters of {IG}
prior, and $|k|$ is the size of the model $k$. By a simple
rearrangement of terms in the above expression, we obtain
\begin{eqnarray*}
&&\hspace*{-6pt}P\bigl(\beta, Z=k \mid Y, \sigma^2\bigr)
\\
&&\hspace*{-6pt}\quad \propto\exp \biggl\{ -\frac{1}{2 \sigma^2} \bigl( (\beta- \tilde{\beta})'
\bigl(D_k + X'X\bigr) (\beta- \tilde{\beta}) - \tilde{
\beta}'\bigl(D_k + X'X\bigr)\tilde {
\beta} \bigr) \biggr\} |D_k|^{{1}/{2}} s_n^{|k|},
\end{eqnarray*}
where $\tilde{\beta} = (D_k + X'X)^{-1} X'Y$. Note that $\tilde
{\beta}$
is a shrinkage estimator of the regression vector $\beta$. Shrinkage of
$\tilde{\beta}$ depends on $D_k$, which is the precision matrix of
$\beta$ given $Z = k$. The components of $\tilde{\beta}_i$
corresponding to $k_i = 0$ are shrunk towards zero while the shrinkage
of coefficients corresponding to $k_i =1$ is negligible (as $\tau
_{1n}^{-2}$ is small).
%
\begin{eqnarray}
\label{post-final} %
P\bigl(Z =k \mid Y,
\sigma^2\bigr) & \propto &Q_k s_n^{|k|}
\exp \biggl\{ -\frac{1}{2 \sigma^2} \bigl(Y'Y - \tilde {
\beta}'\bigl(D_k + X'X\bigr)\tilde{\beta}
\bigr) \biggr\}\nonumber
\\
& =& Q_k s_n^{|k|} \exp \biggl\{ -
\frac{1}{2 \sigma^2} \bigl(Y'Y - Y'X\bigl(D_k
+ X'X\bigr)^{-1}X'Y \bigr) \biggr\}
\\
& =& Q_k s_n^{|k|} \exp \biggl\{ -
\frac{1}{2 \sigma^2} \tilde{R}_k \biggr\}, \nonumber
\end{eqnarray}
where $Q_k = |D_k + X'X|^{-{1}/{2}} |D_k|^{{1}/{2}}$. Next, we
obtain bounds on $Q_k$.

%
\begin{lem}\label{scale-lem} Let $A$ be an invertible matrix, and $B$
be any matrix with appropriate dimension. Further, let $k$ and $j$ be
any pair of models. Then,
\begin{longlist}[(iii)]
\item[(i)]
$ |(A + B'B)^{-1} A| = |I + B A^{-1}B'|^{-1}$,
\item[(ii)]
$ ( I + \tau_{1n}^2 X_k X_k' + \tau_{0n}^2 X_j X_j' )^{-1}
\geq (I + \tau_{1n}^2 X_k X_k' )^{-1} (1 - \xi_n) $,

where $\xi_n = n \tau_{0n}^2 \lambda_M^n = o(1)$, and
\item[(iii)]
$Q_k \leq w'  (n \tau_{1n}^2 \lambda_m^n (1
- \phi_n)  )^ {- ({1}/{2}) (r_{k}^* - r_t)} (\lambda_m^n )^{-
({1}/{2}) |t \wedge k^c|} Q_t$,
where $w' >0$, $r_k = \operatorname{rank}(X_k)$, $r_k^*= r_k \wedge m_n$,
and $\phi_n
= o(1)$.
\end{longlist}
\end{lem}
\begin{pf}
(i) We use the Sylvester's determinant theorem, and the
multiplicative property of the determinant to obtain
\begin{eqnarray*}
\bigl|\bigl(A + B'B\bigr)^{-1}
A\bigr| & = & \bigl|I + A^{-{1}/{2}} B' B A^{-{1}/{2}}\bigr|^{-1}
\\
& = &\bigl |I + B A^{-1} B'\bigr|^{-1}.
\end{eqnarray*}

(ii) By the Sherman--Morrison--Woodbury (SMW) identity,
assuming $A, C$ and $(C^{-1} + D A^{-1} B)$ to be nonsingular,
%
\begin{equation}
\label{smw-eq} (A + BCD)^{-1} = A^{-1} - A^{-1} B
\bigl(C^{-1} + D A^{-1} B\bigr)^{-1} D A^{-1},
\end{equation}
we have, for any vector $a$,
\[
a' \bigl(I + \tau_{1n}^2
X_k X_k' + \tau_{0n}^2
X_{j} X_{j}'\bigr)^{-1} a =
a' G^{-1} a - \tau_{0n}^2 H,
\]
where $G = I + \tau_{1n}^2 X_k X_k' $ and $H = a' G^{-1} X_{j}(I+
\tau
^2_{0n}X_{j}' G^{-1} X_{j})^{-1}X_{j}' G^{-1} a$. Note that
%
\begin{eqnarray}
0 \leq\tau_{0n}^2 H &\leq&
\tau_{0n}^2 a' G^{-1}
X_{j} X_{j}' G^{-1} a
\nonumber
\\[-8pt]
\\[-8pt]
\nonumber
&\leq& n \tau_{0n}^2 \lambda_M^n
a' G^{-1} a,
\end{eqnarray}
where $\lambda_M^n$ is the maximum eigenvalue of the Gram matrix $X'
X/n$. Therefore,
\[
a' \bigl(I + \tau_{1n}^2
X_k X_k' \bigr)^{-1} a \bigl(1 -
n \tau_{0n}^2 \lambda_M^n\bigr)
\leq a' \bigl(I + \tau_{1n}^2 X_k
X_k' + \tau_{0n}^2
X_{j} X_{j}'\bigr)^{-1} a,
\]
and hence (ii) is proved.

(iii)
From part (i) of the lemma, we have
%
\begin{eqnarray}
Q_k &=&\bigl |I + XD_k^{-1}X'\bigr|^{-{1}/{2}}
\nonumber
\\[-8pt]
\\[-8pt]
\nonumber
& = & \bigl|I + \tau_{1n}^2 X_k
X_k' +\tau_{0n}^2
X_{k^c} X_{k^c}' \bigr|^{-
{1}/{2}}.
\end{eqnarray}
Define $A = I + \tau_{1n}^2 X_{k \wedge t} X_{k \wedge t}' +\tau_{0n}^2
X_{{k^c \vee t}^c} X_{{k^c \vee t}^c}' $. Then, by (ii) we have
\[
(1 - \xi_n) \bigl(I + \tau_{1n}^2
X_{k \wedge t} X_{k \wedge t}'^{-1}\bigr) \leq
A^{-1} \leq\bigl(I + \tau_{1n}^2 X_{k \wedge t}
X_{k \wedge t}'\bigr)^{-1}.
\]

This, along with Condition \ref{eigen-cond} implies
\begin{eqnarray*}
\label{ratio-scale1} %
\frac{Q_k}{Q_{k \wedge t}} & = &
\bigl|I + \tau_{1n}^2 X_k X_k'
+\tau_{0n}^2 X_{k^c} X_{k^c}'
\bigr|^{-
{1}/{2}} |A|^{{1}/{2}}
\\
& = & \bigl|A + \bigl(\tau_{1n}^2 - \tau_{0n}^2
\bigr) X_{k \wedge t^c} X_{k \wedge
t^c}'\bigr|^{-{1}/{2}}
|A|^{{1}/{2}}
\\
& = & \bigl|I + \bigl(\tau_{1n}^2 - \tau_{0n}^2
\bigr) X_{k \wedge t^c} ' A^{-1} X_{k
\wedge t^c}
\bigr|^{- {1}/{2}}
\\
& \leq& \bigl|I + \bigl(\tau_{1n}^2 - \tau_{0n}^2
\bigr) (1 - \xi_n) X_{k \wedge t^c} ' \bigl(I +
\tau_{1n}^{2} X_{k \wedge t} X_{k \wedge t}'
\bigr)^{-1} X_{k \wedge
t^c}\bigr |^{- {1}/{2}}
\\
& = & \bigl|I + \tau_{1n}^{2} X_t
X_t' + \tau_{1n}^2 (1 -
\phi_n) X_{k
\wedge t^c} X_{k \wedge t^c}'\bigr|^{-{1}/{2}}
\bigl|I + \tau_{1n}^{2} X_{k
\wedge t} X_{k \wedge t}'
\bigr|^{{1}/{2}}
\\
& \leq&\bigl |I + \tau_{1n}^2 (1 - \phi_n)
X_k X_k'\bigr|^{-{1}/{2}}\bigl|I +
\tau_{1n}^2 X_{k \wedge t} X_{k \wedge t}'\bigr|^{{1}/{2}}
\\
&\leq& \bigl(n \tau_{1n}^2 \lambda_m^n
(1 - \phi_n) \bigr)^ {- (r_k^* - r_{t
\wedge k})/2} (1 - \phi_n)^{-|t \wedge k|/2},
\end{eqnarray*}
where $(1 - \phi_n) = (\tau_{1n}^2 - \tau_{0n}^2) (1 - \xi_n) /
\tau
_{1n}^2 \rightarrow1$.
Similarly, let $A = I + \tau_{1n}^2 X_{t} X_{t}' +\tau_{0n}^2 X_{t^c}
X_{t^c}'$ to obtain
\begin{eqnarray*}
\label{ratio-scale2} %
\frac{Q_{k \wedge t}}{Q_t}
& = & \bigl|A - \bigl(\tau_{1n}^2 - \tau_{0n}^2
\bigr) X_{k \wedge t^c} X_{k \wedge
t^c}'\bigr|^{-{1}/{2}}
|A|^{{1}/{2}}
\\
& \leq& \bigl|I + \tau_{1n}^{2} X_{k \wedge t}
X_{k \wedge t}'\bigr|^{-{1}/{2}} \bigl|I + \tau_{1n}^{2}
X_t X_t' \bigr|^{{1}/{2}}
\\
&\leq&\bigl |I + \tau_{1n}^{2} X_{t \wedge k^c}
X_{t \wedge k^c}'\bigr|^{{1}/{2}}
\\
& \leq& \bigl(n \tau_{1n}^2 c'
\bigr)^{|t \wedge k^c|/2}.
\end{eqnarray*}
The above two inequalities give
\[
\label{ratio-scale} %
 \frac{Q_k}{Q_t}
\leq w' \bigl(n \tau_{1n}^2
\lambda_m^n (1 - \phi_n)
\bigr)^ {- (r_k^* -
r_t)/2} \bigl({\lambda_m^n}
\bigr)^{-|t \wedge k^c|/2}.
\]
\upqed\end{pf}
Due to \eqref{post-final}, we have
\[
\operatorname{PR}(k, t)  = \frac{Q_k}{Q_t}
s_n^{|k| -|t|} \exp \biggl\{ -\frac{1}{2 \sigma
^2} (
\tilde{R}_k - \tilde{R}_t) \biggr\} .
\]
Therefore, Lemma~\ref{scale-lem}(iii) implies Lemma~\ref{bf-lem}.
\end{pf*}

\begin{pf*}{Proof of Lemma~\ref{mnev-lemma}}
The rows of $X_k$ are $n$ independent sub-Gaussian random isotropic
random vectors in $R^{|k|}$. Note that $|k| \leq m_n$ implies $|k| =
o(n)$. Due to Theorem~5.39 of \citet{Vershynin12}, with probability at
least $ 1 - 2 \exp(-c s)$, we have
%
\begin{equation}
\label{mineig-eq} %
\phi_{\mathrm{min}} \biggl(
\frac{X_k' X_k}{n} \biggr) > \biggl(1 - C \sqrt {\frac
{|k|}{n}} - \sqrt{
\frac{s}{n}} \biggr)^2,
\end{equation}
where $c$ and $C$ are absolute constants that depend only on the
sub-Gaussian norms of the rows of the matrix $X_k$.

Let us fix $s= n(1 - \phi)$ for some $\phi>0$, and define the event
given by equation~\eqref{mineig-eq} as $A_k$. We then have $P[A_k^c] <
2 \exp(-c (1 - \phi) n)$ for all $k$. By taking an union bound over
$\{
k\dvtx |k| \leq m_n \}$, we obtain
\begin{eqnarray*}
 P\biggl[\bigcup
_{|k| \leq m_n} A_k^c\biggr] &\leq&
p_n^{ m_n} \exp\bigl(-c (1 - \phi) n\bigr)
\\
&=& \exp{ \biggl\{ \frac{n}{2 + \nu} -c (1 - \phi) n \biggr\}} \rightarrow0,
\end{eqnarray*}
if $ \nu> (\frac{1}{c (1 - \phi)} - 2)$.
Therefore, in the event $\bigcap_{|k| \leq m_n} A_k$, whose probability
goes to~1, we have $\phi_{\mathrm{min}}  (X_k' X_k/n ) \geq\phi
^2/4 -
O(\sqrt{m_n/n})>0$, for all $k$.
\end{pf*}

\section*{Acknowledgments}
The authors are grateful to anonymous
referees and an Associate Editor for their encouraging and helpful
comments on an earlier version of the paper. The authors would also
like to thank Professors Howard Bondell, Val Johnson and Faming Liang
for sharing with us their code to perform Bayesian model selection.

\begin{supplement}[id=suppA]
\stitle{Supplement to ``Bayesian variable selection with shrinking and
diffusing priors''}
\slink[doi]{10.1214/14-AOS1207SUPP} 
\sdatatype{.pdf}
\sfilename{aos1207\_supp.pdf}
\sdescription{This supplement contains the proofs of Theorems \ref
{sigfix-thmm}, \ref{sigprior-thmm} and Lemma~\ref{rates-lem}.}
\end{supplement}

%

%

\printaddresses

\end{document}